\documentclass{birkjour}
\usepackage{amsmath, amsfonts, amsthm, amssymb, setspace, textcomp, bbm, multirow, tikz,mathtools,enumitem,booktabs,longtable}

\newtheorem{theorem}{Theorem}
\numberwithin{theorem}{section}
\numberwithin{equation}{section}

\theoremstyle{definition}

\theoremstyle{remark}
\newtheorem*{remark}{Remark}

\newcommand{\N}{\mathbb N}
\newcommand{\Z}{\mathbb Z}

\begin{document}
\allowdisplaybreaks
%
%
%
%
%
%
%
%
%

\title{A Note on Andrews' Partitions with Parts Separated by Parity}

\author[Bringmann]{Kathrin Bringmann}
\address{Mathematical Institute, University of Cologne, Weyertal 86-90, 50931 Cologne, Germany}
\email{kbringma@math.uni-koeln.de}
\thanks{The research of the first author is supported by the Alfried Krupp Prize for Young University Teachers of the Krupp foundation and the research leading to these results receives funding from the European Research Council under the European Union's Seventh Framework Programme (FP/2007-2013) / ERC Grant agreement n. 335220 - AQSER}

\author[Jennings-Shaffer]{Chris Jennings-Shaffer}
\address{Mathematical Institute, University of Cologne, Weyertal 86-90, 50931 Cologne, Germany}
\email{cjenning@math.uni-koeln.de}



\dedicatory{Dedicated to George Andrews in honor of his 80th birthday}

\begin{abstract}
In this note we give three identities for partitions with parts 
separated by parity, which were recently introduced by Andrews.
\end{abstract}

\maketitle

\section{Introduction}

Recently Andrews \cite{Andrews1} studied integer partitions in which all parts of a given parity
are smaller than those of the opposite parity. Furthermore, he considered
eight subcases based on the parity of the smaller parts
and parts of a given parity appearing at most once or an unlimited number of times.
Following Andrews, we use ``ed'' for evens distinct, ``eu'' for evens unlimited,
``od'' for odds distinct, and ``ou'' for odds unlimited. With
``zw'' and ``xy'' from the four choices above, we let
$F^{\rm zw}_{\rm xy}(q)$ denote the generating function of partitions where
zw specifies the parity and condition of the larger parts and
xy specifies the parity and condition of the smaller parts.

The eight relevant generating functions are
\begin{align*}
F^{\rm ou}_{\rm eu}(q) 
&:=
	\sum_{n=0}^\infty \frac{q^{2n}}{\left(q^2;q^2\right)_n\left(q^{2n+1};q^2\right)_\infty}
,\\
F^{\rm od}_{\rm eu}(q) 
&:=
	\sum_{n=0}^\infty \frac{q^{2n} \left(-q^{2n+1};q^2\right)_\infty }{\left(q^{2};q^2\right)_n}
,\\
F^{\rm ou}_{\rm ed}(q) 
&:=
	\sum_{n=0}^\infty \frac{\left(-q^2;q^2\right)_nq^{2n+2}}{\left(q^{2n+3};q^2\right)_\infty}
,\\
F^{\rm od}_{\rm ed}(q) 
&:=
	\sum_{n=0}^\infty q^{2n+2}\left(-q^2,q^2\right)_n\left(-q^{2n+3};q^2\right)_\infty
,\\
F^{\rm eu}_{\rm ou}(q) 
&:=
	\sum_{n=0}^\infty \frac{q^{2n+1}}{\left(q;q^2\right)_{n+1}\left(q^{2n+2};q^2\right)_\infty}
,\\
F^{\rm ed}_{\rm ou}(q)
&:= 
	\sum_{n=0}^\infty \frac{q^{2n+1}\left(-q^{2n+2};q^2\right)_\infty}{\left(q;q^2\right)_{n+1}}
,\\
F^{\rm eu}_{\rm od}(q) 
&:=
	\sum_{n=0}^\infty \frac{q^{2n+1}\left(-q;q^2\right)_n}{\left(q^{2n+2};q^2\right)_\infty}
,\\
F^{\rm ed}_{\rm od}(q) 
&:=
	\sum_{n=0}^\infty q^{2n+1}\left(-q;q^2\right)_n\left(-q^{2n+2};q^2\right)_{\infty}
.
\end{align*}
Here we are using the standard product notation 
$(a;q)_n := \prod_{j=0}^{n-1}(1-aq^j)$ for $n\in\N_0\cup\{\infty\}$.
We note that with the exception of $F^{\rm ou}_{\rm eu}(q)$ and 
$F^{\rm od}_{\rm eu}(q)$, we do not allow the subpartition consisting
of the smaller parts to be empty.

Andrews' identities (after minor corrections) can be stated as
\begin{align*}
F^{\rm ou}_{\rm eu}(q) 
&=
	\frac{1}{\left(1-q\right)\left(q^2;q^2\right)_\infty}
,\\
F^{\rm od}_{\rm eu}(q) 
&=
	\frac{1}{2}\left(\frac{1}{\left(q^2;q^2\right)_\infty}+\left(-q;q^2\right)_\infty^2 \right)
,\\
F^{\rm ou}_{\rm ed}(-q) 
&=
	\frac{1}{2\left(-q;q^2\right)_\infty}\left(
		\left(-q;q\right)_\infty - 1 - \sum_{n=0}^\infty q^{\frac{n\left(3n-1\right)}{2}}\left(1-q^n\right)
	\right)
,\\
F^{\rm eu}_{\rm ou}(q) 
&=
	\frac{1}{1-q}\left(
		\frac{1}{\left(q;q^2\right)_\infty} - \frac{1}{\left(q^2;q^2\right)_\infty}
	\right)
,\\
F^{\rm ed}_{\rm ou}(-q)
&= 
	-\frac{\left(-q^2;q^2\right)_\infty}{2}\left(
		2 - \frac{1}{\left(-q;q\right)_\infty} 
		-
		\sum_{n=0}^\infty \frac{q^{n^2+n}}{\left(-q;q\right)_n^2\left(1+q^{n+1}\right)}		
	\right)
,\\
F^{\rm eu}_{\rm od}(-q) 
&=
	-\frac{1}{\left(q^2;q^2\right)_\infty}\sum_{j=1}^\infty\sum_{n=j}^\infty
	\left(-1\right)^{n+j}q^{\frac{n\left(3n+1\right)}{2}-j^2}\left(1-q^{2n+1}\right)
.
\end{align*}
Surprisingly, these identities are derived with little more than the $q$-binomial theorem,
Heine's transformation, and the Rogers-Fine identity. In the following
theorem, we give new identities for $F^{\rm od}_{\rm ed}(q)$, 
$F^{\rm ed}_{\rm od}(q)$, and $F^{\rm ed}_{\rm ou}(-q)$.

\begin{theorem}\label{TheTheorem}
The following identities hold,
\begin{align}
\label{Eqoded}
F^{\rm od}_{\rm ed}(q) 
&=
	\frac{q\left(-q;q^2\right)_\infty}{1-q}\left(
		1
		-
		\frac{(-q^2;q^2)_\infty}{(-q;q^2)_\infty}
	\right)
,\\
\label{Eqedod}
F^{\rm ed}_{\rm od}(q) 
&=
\frac{q(-q^2;q^2)_\infty}{1-q}\left( 2 - \frac{(-q;q^2)_\infty}{(-q^2;q^2)_\infty}\right)
,\\
\label{Eqedou}
F^{\rm ed}_{\rm ou}(-q)
&= 
	-\frac{\left(-q^2;q^2\right)_\infty}{2}\left(
		2 - \frac{1}{\left(-q;q\right)_\infty}
		-\frac{2}{\left(q;q\right)_\infty}\sum_{n\in\Z}
		\frac{\left(-1\right)^nq^{\frac{3n\left(n+1\right)}{2}}}{1+q^n}
	\right)
.
\end{align}
\end{theorem}
\begin{remark}
The functions $F^{\rm od}_{\rm ed}(q)$ and $F^{\rm ed}_{\rm od}(q)$ are basically
modular functions. Also we find that $F^{\rm ed}_{\rm ou}(-q)$ is related 
to Ramanujan's third order mock theta function $f(q)$, as
\begin{align*}
f(q) &:=
\sum_{n=0}^\infty \frac{q^{n^2}}{(-q;q)_n^2}
= 
\frac{2}{\left(q;q\right)_\infty}\sum_{n\in\Z}
	\frac{\left(-1\right)^nq^{\frac{n\left(3n+1\right)}{2}}}{1+q^n}
\\&=
2-\frac{2}{\left(q;q\right)_\infty}\sum_{n\in\Z}
	\frac{\left(-1\right)^nq^{\frac{3n\left(n+1\right)}{2}}}{1+q^n}
,
\end{align*}
where the final equality uses Euler's pentagonal numbers theorem.
\end{remark}

\section{Proof of Theorem \ref{TheTheorem}}

To prove equations \eqref{Eqoded} and \eqref{Eqedod}, we require the following
$q$-series identity,
\begin{align}\label{EqQseries1}
\sum_{n=0}^\infty \frac{(x;q)_nq^n}{(y;q)_n}
&=
\frac{q(x;q)_\infty}{y(y;q)_\infty \left(1-\frac{xq}{y}\right)}
+
\frac{\left(1-\frac{q}{y}\right)}{\left(1-\frac{xq}{y}\right)}.
\end{align}
We note that \eqref{EqQseries1} is (4.1) from \cite{AndrewsSubbaraoVidyasagar1}
and was proved with  
Heine's transformation
\cite[page 241, (III.2)]{GasperRahman1}.
To prove equation \eqref{Eqedou} we require the
concept of a Bailey pair and Bailey's Lemma, which are described in
\cite[Chapter 3]{Andrews3}. A pair of sequences
$(\alpha,\beta)$ is called a \textit{Bailey pair relative} to $a=q$ if 
\begin{align*}
\beta_n &= \sum_{j=0}^n \frac{\alpha_j}{(q;q)_{n-j}(q^2;q)_{n+j}}.
\end{align*}
A limiting form of Bailey's Lemma states that if $(\alpha_n,\beta_n)$
is a Bailey pair relative to $q$, then
\begin{align}\label{EqBaileysLemma}
\sum_{n=0}^\infty q^{n^2+n}\beta_n
&=
\frac{1}{(q^2;q)_\infty}\sum_{n=0}^\infty q^{n^2+n}\alpha_n.
\end{align}
The Bailey pair we use is given by
\begin{align}\label{EqBaileyPair}
\beta_n^\prime &:= \frac{1}{(-q;q)_n^2(1+q^{n+1})}
,&
\alpha_n^\prime &:= \frac{2(-1)^n  q^{\frac{n(n+1)}{2}}(1-q^{2n+1})}
	{(1-q)(1+q^n)(1+q^{n+1})}
,
\end{align}
which follows from taking the Bailey pair
from Theorem 8 of \cite{Lovejoy1} with
$a\rightarrow q$, $b=-1$, $c=-q$, and $d=-1$
and dividing both $\alpha_n$ and $\beta_n$
by $(1+q)$.

\begin{proof}[Proof of Theorem \ref{TheTheorem}]
We find that
\begin{align*}
F^{\rm od}_{\rm ed}(q)
&=
\left(-q;q^2\right)_\infty\sum_{n=1}^\infty \frac{\left(-q^2;q^2\right)_{n-1}q^{2n}}{\left(-q;q^2\right)_n}
\\
&=
\frac{\left(-q;q^2\right)_\infty}{2}\left(
	-1
	+\sum_{n=0}^\infty \frac{\left(-1;q^2\right)_{n}q^{2n}}{\left(-q;q^2\right)_n}
\right).
\end{align*}
With $q\mapsto q^2$, $x=-1$, and $y=-q$, equation \eqref{EqQseries1} implies that
\begin{align*}
\sum_{n=0}^\infty \frac{(-1;q^2)q^{2n}}{(-q;q^2)_n}
&=
-\frac{q\left(-1;q^2\right)_\infty}{\left(-q;q^2\right)_\infty(1-q)} + \frac{1+q}{1-q}.
\end{align*}
Equation \eqref{Eqoded} then follows after elementary simplifications.

Similarly, we have that
\begin{align*}
F^{\rm ed}_{\rm od}(q)
&=
\left(-q^2;q^2\right)_\infty\sum_{n=0}^\infty \frac{\left(-q;q^2\right)_{n}q^{2n+1}}{\left(-q^2;q^2\right)_n}.
\end{align*}
By applying \eqref{EqQseries1} with $q\mapsto q^2$, $x=-q$, and $y=-q^2$, we find that
\begin{align*}
\sum_{n=0}^\infty \frac{\left(-q;q^2\right)q^{2n}}{\left(-q^2;q^2\right)_n}
&=
-\frac{\left(-q;q^2\right)_\infty}{\left(-q^2;q^2\right)_\infty(1-q)} + \frac{2}{1-q},
\end{align*}
and \eqref{Eqedod} follows.

For $F^{\rm ed}_{\rm ou}(q)$, we begin with Andrews' identity \cite{Andrews1}
\begin{align*}
F^{\rm ed}_{\rm ou}(-q)
&=
-\frac{\left(-q^2;q^2\right)_\infty}{2}\left(
	2- \frac{1}{(-q;q)_\infty}	
	-\sum_{n=0}^\infty \frac{q^{n^2+n}}{(-q;q)_n^2\left(1+q^{n+1}\right)}
\right).
\end{align*}
By applying \eqref{EqBaileysLemma} to the Bailey pair 
$(\alpha^\prime,\beta^\prime)$ in \eqref{EqBaileyPair},
we have that
\begin{align*}
\sum_{n=0}^\infty \frac{q^{n^2+n}}{(-q;q)_n^2\left(1+q^{n+1}\right)}
&=
\frac{2}{(q;q)_\infty}
\sum_{n=0}^\infty \frac{(-1)^nq^{\frac{3n(n+1)}{2}} \left(1-q^{2n+1}\right)}{\left(1+q^n\right)\left(1+q^{n+1}\right)}.
\end{align*}
We use the partial fraction decomposition
\begin{align*}
\frac{1-q^{2n+1}}{\left(1+q^n\right)\left(1+q^{n+1}\right)}
&=
\frac{1}{1+q^n}-\frac{q^{n+1}}{1+q^{n+1}}
,
\end{align*}
to deduce that
\begin{align*}
\sum_{n=0}^\infty \frac{(-1)^nq^{\frac{3n(n+1)}{2}} \left(1-q^{2n+1}\right)}{\left(1+q^n\right)\left(1+q^{n+1}\right)}
&=
\sum_{n=0}^\infty (-1)^nq^{\frac{3n(n+1)}{2}}
\left( \frac{1}{1+q^n}-\frac{q^{n+1}}{1+q^{n+1}}\right)
\\
&=
\sum_{n\in\Z} \frac{(-1)^nq^{\frac{3n(n+1)}{2}}}{1+q^n}.
\end{align*}
Altogether this implies equation \eqref{Eqedou}.
\end{proof}

By applying Theorem 1.1 part $3$ of \cite{Lovejoy2} to the Bailey
pair $E(3)$ of \cite{Slater1}, 
we find that
\begin{multline*}
F^{\rm ed}_{\rm od}(-q) =
-\frac{q\left(q;q\right)_\infty\left(-q^2;q^2\right)_\infty}{\left(q^2;q^2\right)_\infty^2}
\\
\times \sum_{n=0}^\infty\sum_{m=0}^\infty \left(-1\right)^mq^{\frac{n\left(n+3\right)}{2}+2nm+2m^2+2m}\left(1+q^{2m+1}\right)
.
\end{multline*}
As such, we have that
\begin{multline*}
\left(\sum_{n,m\geq 0} - \sum_{n,m<0}\right) (-1)^m q^{\frac{n(n+3)}{2}+2nm +2m(m+1)}\\
=
\frac{2\left(q^2;q^2\right)_\infty}{(1+q)\left(q;q^2\right)_\infty}-\frac{\left(q^2;q^2\right)_\infty}{(1+q)\left(-q^2;q^2\right)_\infty}.
\end{multline*}
We note that the corresponding quadratic form is degenerate, and so a priori the
modularity properties of this theta function are unclear. More generally, one
can prove directly that, for $c\in\N$,  
\begin{align*}
\sum_{n,m\geq 0}z^nw^mq^{n^2+2cnm +c^2m^2}
&=
\frac{1}{1-\frac{w}{z^c}}\sum_{k=0}^{c-1}\sum_{n=0}^\infty
z^{cn+k}q^{(cn+k)^2}
\left(1-\frac{w^{n+1}}{z^{cn+c}}\right).
\end{align*}
The above is a sum of partial theta functions, which sometimes
combine to give a modular form.

\section*{Acknowledgments}
The authors thank George Andrews, Karl Mahlburg, and the anonymous referee for their careful reading
and comments on an earlier version of this manuscript.

\end{document}